\input amstex
\magnification=\magstep1
\input epsf
\input amssym.def
\input amssym
\pageno=1
\baselineskip 14 pt
\def \pop#1{\vskip#1 \baselineskip}

\font\gr=cmbx12
\def \nr{\operatorname {nr}}
\def \sep{\operatorname {sep}}

\def \et{\operatorname {et}}
\def \dim{\operatorname {dim}}

\def \Spec{\operatorname {Spec}}
\def \Spf{\operatorname {Spf}}
\def \lim{\operatorname {lim}}

\def \Sp  {\operatorname {Sp  }}

\def \fppf{\operatorname {fppf }}
\def \Pic{\operatorname {Pic }}
\def \ord{\operatorname {ord }}
\def \ord{\operatorname {ord }}

\def \res{\operatorname {res }}

\def \log{\operatorname {log}}

\def \cond{\operatorname {cond}}

\def \Spf{\operatorname {Spf}}

\def \Aut{\operatorname {Aut}}
\def \Hom{\operatorname {Hom}}

\def \kum{\operatorname {kum}}
\def \Ker{\operatorname {Ker}}
\def \Fil{\operatorname {Fil}}

\def \SS-Deg{\operatorname {SS-Deg}}
\def \SNS-Deg{\operatorname {SNS-Deg}}
 \def \DS-Deg{\operatorname {DS-Deg}}
\def \DNS-Deg{\operatorname {DNS-Deg}}

\pop {4}
\par
\noindent                                          
\centerline {\bf \gr Torsors under finite and flat group schemes of rank $p$} 
\par
\noindent                                          
\centerline {\bf \gr with Galois action}

\pop {3}
\noindent                                          
\centerline {\bf \gr Mohamed Sa\"\i di}

\pop {4}
\centerline {\bf \gr Abstract}
\pop {.5}
\par
In this note we study the geometry of torsors under flat and finite 
commutative group schemes
of rank $p$ above curves in characteristic $p$, and above 
relative curves over a complete 
discrete valuation ring of inequal characteristics. 
In both cases we study the Galois action of the Galois group of 
the base field on these torsors. We also study the
degeneration of $\mu_p$-torsors from characteritic $0$ to characteristic $p$
and show that this degeneration is compatible with the Galois action.
We then discuss the lifting of torsors under flat and commutative group schemes
of rank $p$ from positive to zero characteristics. Finally, for a proper
and smooth curve $X$ over a complete discrete valuation field of inequal 
characteristics we show the existence of a canonical Galois equivariant 
filtration on the first \'etale cohomology group of the geometric fibre of $X$
with values in $\mu_p$.
\pop {2}
\par
\noindent
{\bf \gr 0. Introduction.}\rm \
In this note we study the geometry of torsors under flat and finite 
commutative group schemes of rank $p$ above curves in 
characteristic $p$ (cf. I), and above curves over a discrete complete 
valuation field of inequal characteristics (cf. II). For a connected scheme
$Z$ over a field of characteristic $p>0$ we define the group
$H^1_{\fppf}(Z)_p$ of ``mixed torsors'' of rank $p$, and in case
$Z$ is a semi-stable curve we define the group
$H^1_{\fppf}(Z)_p^{\kum}$ which 
classify ``mixed kummerian torsors'' of rank $p$. These torsors 
arise naturally in the case of degeneration of $\mu_p$-torsors 
from characteristic $0$ to characteristic $p$ (cf. 1.7 and [Sa-1]). In all 
these cases we study the canonical action of the Galois group of 
the base field on these torsors.

\par
in paragraph II we consider a discrete complete valuation ring $R$ of inequal 
characteristics with fraction field $K$ and residue field $k$, 
and a formal $R$-scheme $X$ of finite type which is normal 
flat over $R$, and smooth of relative dimension 1. We denote by 
$G_K$ the Galois group of a separable
closure of $K$. Let $X_{\overline K}$ (resp. $X_{\overline k}$) denote the 
geometric generic fibre of $X$ viewed as a rigid analytic space 
(resp. the geometric special fibre of $X$). We study the degeneration of 
$\mu_p$-torsors of $X_{\overline K}$. Our main result is the following
which in particular states that the degeneration of $\mu_p$-torsosrs of
$X_{\overline K}$ is compatible with the 
Galois action:

\pop {.5}
\par
\noindent                                          
{\bf \gr Theorem.(cf. 2.5.1 and 2.7)}\rm \ {\sl There exists a canonical
specialisation  group homomorphism:
$$\Sp: H^1_{\et}(X_{\overline K},\mu_p)\to 
H^1_{\fppf}(X_{\overline k},\mu_p)$$ 
which is $G_K$-equivariant. Moreover there exists a canonical
specialisation map (this is not a group homomorphism):
$$\Sp: H^1_{\et}(X_{\overline K},\mu_p)\to 
H^1_{\fppf}(X_{\overline k})_p$$ which is $G_K$-equivariant.}
\pop {.5}
\par
Hier we consider the groups $H^1_{\fppf}(X_{\overline k},\mu_p)$
and $H^1_{\fppf}(X_{\overline k})_p$ as $G_K$-sets via the canonical quotient
$G_K\to G_k$ of $G_K$, where $G_k$ is the Galois group of the residue field 
$k$ of $R$, and via the canonical action of $G_k$ on these groups. Also
by construction the kernel of $\Sp: H^1_{\et}(X_{\overline K},\mu_p)\to 
H^1_{\fppf}(X_{\overline k},\mu_p)$ corresponds to those $\mu_p$-torsors 
of $X_{\overline K}$
which degenerate to either an $\alpha_p$ or an \'etale torsor of
$X_{\overline k}$.

\pop {.5}
\par
in paragraph III we explain the degeneration of $\mu_p$-torsors 
above the boundaries 
of formal fibres of formal $R$-curves at closed points, and define their 
degeneration type. In paragraph IV we discuss the lifting of torsors 
under flat and 
finite commutative group schemes of rank $p$ from characteristic $p$
to characteristic $0$, and show in particular that such torsors above
proper curves can alaways be lifted to characteristic zero (cf. 4.5). 
Finally for a scheme $X$ which is proper, smooth, of relative dimension 1
over a discrete complete valuation ring of inequal characteristics
we define a canonical decreasing filtration on the group 
$H^1_{\et}(X_{\bar \eta},\mu_p)$, where $X_{\bar \eta}$ is the geometric 
generic fibre of $X$ (cf. 5.1). Moreover we show that this filtration is 
equivariant under the action of the Galois group of the base field (cf. 5.2).

\pop {.5}
\par
This paper is the first of a serie of papers [Sa] and [Sa-1] where we 
compute the vanishing cycles arising from the degeneration of 
$\mu_p$-torsors above curves as well as the semi-stable reduction of 
these torsors and the Galois action on them.

\pop {.5}
\par
\noindent
{\bf \gr I. Torsors under finite and flat group schemes of rank $p$ in 
characteristic $p$.}

\pop {.5}
\par
In this section we discuss the geometry of torsors under finite and flat 
group schemes of rank $p$ above curves in characteristic $p$. For such
torsors we introduce the notion of ``conductor'' and ``residue'' at
the ``critical points'', and we carefully explain 
the Galois action on them and the fact that conductor and residue are 
invariant under this action. We also introduce for a semi-stable curve
the group $H^1_{fppf}(\ )^{\kum}$ which classifies 
``{\it kummerian mixed torsors}'' and which appear naturally when considering
the degeneration of $\mu_p$-torsors from zero to positive characteristics
(cf. [Sa-1]).

\pop {.5}
\par
\noindent
{\bf \gr 1.1. Kummer theory.}\rm
\par
Let $X$ be a scheme and let $l$ be a positif
integer. The following sequence
is exact on $X_{\fppf}$:
$$(1)\ \ 1\to \mu _l\to \Bbb G_m @>{x^l}>> \Bbb G_m\to 1$$
\par
Let $f:Y\to X$ be a $\mu_l$-torsor. Then there exists an open covering 
$\Cal U:=(U_s)_s$ of $X$ and 
invertible fonctions $u_s\in \Gamma (U_s,\Cal O_X)^*$, which are 
defined up to multiplication by
a $p$-power, such that above $U_s$ the torsor $f$ is 
given by an equation 
$T_s^p=u_s$ (cf. [Mi], III, 4).

\pop {.5}
\par
Assume now that $X$ is a scheme over a field $k$. Let 
$\overline k$ be an algebraic closure of $k$, and let
$G_k$ be the Galois group of the separable closure $k^{\sep}$ of $k$ in 
$\overline k$ of $k$. The 
group $G_k$ acts canonically by automorphisms 
on $\overline X:=X\times _k{\overline k}$, and this action induces 
in a natural way an action of $G_k$ on $\mu_l$-torsors of $\overline X$. 
There exists a canonical homomorphism:
$$(1')\ \ G_k\to \Aut H_{\fppf}^1(\overline X,\mu_l)$$
\par
More precisely, with the same notations as above,  an element 
$\sigma\in G_k$ acts on $\overline X$ hence on 
coverings of $\overline X$. The element $\sigma$ associates to the covering 
$\Cal U:=(U_s)_s$ the covering $\Cal U ^{\sigma}:={(U^{\sigma}_s)}_s$. To a 
$\mu_l$-torsor $f:\overline Y\to \overline X$ is then associated 
the $\mu_l$-torsor $f^{\sigma}:\overline Y^{\sigma}\to X$,
which is locally defined by an equation $T_s^p=u^{\sigma}_s$, where
$u^{\sigma}_s$ is the image via $\sigma$ of $u_s$ which is a unit on 
$U^{\sigma}_s$.

\pop {1}
\par
\noindent
{\bf \gr 1.2.}\rm\  In what follows and unless otherwise is specified, we 
will use the following
notations: let $U$ be a smooth geometrically connected algebraic curve 
over a field $k$ of characteristic $p>0$, and let $X$ be the smooth
compactification of $U$. Let $G_k$ denote the
Galois group of a separable closure $k^{\sep}$ of $k$ in a fixed algebraic 
closure $\overline k$ of $k$. Let $\overline U:=U\times _k{\overline k}$,
and denote by $\overline X$ the smooth compactification of $\overline U$. 
If $U\neq X$, let $S:=X-U$ and let
$\overline S:=S\times _k{\overline k}$. 

\pop {.5}
\par
\noindent
{\bf\gr 1.3. $\Bbb Z/p \Bbb Z$-Torsors in characterisitc $p>0$.}\rm\ 
The following sequence is exact for the \'etale topology:
$$(2)\ \ 0\to \Bbb Z/p \Bbb Z \to \Bbb G_a@>{x^p-x}>>
\Bbb G_a\to 0$$
\par
Let $f:V\to U$ be a non trivial $\Bbb Z/p \Bbb Z$-torsor. Let $Y$ 
be the smooth compactification of $V$ and let $f':Y\to X$ be the 
canonical morphism,
which is finite of degree $p$, generically separable, and eventually 
ramified above $X-U$. There exists
an open covering $(U_s)_s$ of $U$ and regular fonctions 
$a_s\in \Gamma (U_s,\Cal O_X)$, which are defined up to addition of elements 
of the form $b^p-b$, such that above $U_s$ the torsor $f$ is 
given by an equation $T_s^p-T_s=a_s$. Above a point $x_i\in S$, and after an 
\'etale localisation, the morphism $f'$ is given by an equation $T_i^p-T_i=
t_i^{m_i}$ 
where $m_i$ is an integer prime to $p$, and $t_i$ is a uniformising 
parameter at $x_i$. If $m_i$ is positif then $f'$ is 
\'etale above $x_i$. Otherwise $f'$ ramify above $x_i$, and
the contribution to the degree of the different of the point $x_i$ is 
$(-m_i+1)(p-1)$. In the later case we call the integer 
$\cond _{x_i}(f):=-m_i$ the
{\bf conductor} of the above torsor $f$ at the point $x_i$ (this is the 
classical Hasse conductor). In the case 
where $f$ is \'etale
above $x_i$ we define $\cond _{x_i}(f):=0$ to be the {\bf conductor} 
of $f$ at $x_i$. We define the {\bf residue} $h_i:=\res _{x_i}(f)$ of 
$f$ at the point $x_i$ in any case to be $0$. 
\pop {.5}
\par
The Galois group $G_k$ acts in a natural way on 
$\Bbb Z/p \Bbb Z$-torsors of $\overline U$. More precisely there exists a
canonical homomorphism:
$$(2')\ \ G_k\to \Aut H_{\et}^1(\overline U,\Bbb Z/p \Bbb Z)$$
\par
Let $\sigma \in G_k$. Let $f:\overline V\to \overline U$ be an \'etale 
$\Bbb Z/p\Bbb Z$-torsor, and let $f^{\sigma}:\overline V^{\sigma}\to 
\overline U$ be the $\Bbb Z/p\Bbb Z$-torsor associated to $f$ via the above 
action $(2')$. If $f$ is locally given by an equation $T_s^p-T_s=a_s$, 
then $f^{\sigma}$ is locally given by the equation $T_s^p-T_s=a_s^{\sigma}$,
where $a_s^{\sigma}$ is the image of $a_s$ via $\sigma$ which 
is a regular function on $U_s^{\sigma}$.
Moreover for a point $x_i\in \overline S$ we have 
$\cond _{x_i}(f)=\cond _{x_i^{\sigma}}(f^{\sigma})$ (resp. 
$\res _{x_i}(f)=\res _{x_i^{\sigma}}(f^{\sigma})$), where
$x_i^{\sigma}\in \overline S$ is the image of $x_i$ via $\sigma$.

\pop {.5}
\noindent
{\bf \gr 1.4. $\mu _p$-Torsors.} 
\rm
Let $f:V\to U$ be a $\mu_p$-torsor. Then $f$ is a finite radicial morphism 
of degree $p$. In particular $V$ is homeomorphic to $U$. Let 
$(U_s)_s$ and $u_s\in \Gamma (U_s,\Cal O_X)^*$ be as in 2.1.
The logarithmic differential form $\omega:=du_s/u_s$ of $u_s$ is a
global differential form on $X$, it is the 
{\bf differential form} associated 
to the torsor $f$ (it is well defined). Let 
$\{y_j\}_{j\in Z}\subset U$ be 
the set of zeros of $w$ 
which are contained in $U$, and let $m_j-1:=\ord _{y_j}(\omega)$, with
$m_j\ge 2$ an integer prime to $p$. Locally for the \'etale topology 
the torsor $f$ is defined above the point $y_j$ by an equation 
$T^p=t_j^{m_{j}}$ where $t_j$ is a uniformising 
parameter at $y_j$. The unique point 
$z_j\in V$ above the point $y_j$ is a singular point of $V$, and $V$ is 
unibranche at $y$. Moreover the only singularities of $V$ lie above the 
zeros of $w$ which are contained in $U$. For a point $x\in X$ we 
call the integer $\cond _{x}(f):=m_x:=-(\ord_x(\omega)+1)$ the {\bf conductor} 
of the torsor $f$ at the point $x$. We define the {\bf residue} 
$h_x\in \Bbb F_p$ of the torsor $f$ at the point $x$ to be $0$, unless $m_x=0$,
in which case $h_x:=\res _x(\omega)$. The differential form $\omega$
is fixed under the Cartier operation $C$. In fact the following sequence
is exact on $U_{\et}$:
$$(3)\ 0\to \Cal O_U^{*}@>{F}>> \Cal O_U^{*}@>{d\log}>>
\Omega ^1_{X}@>{C-1}>> \Omega ^1_X\to 0$$
where $F$ is the Frobenius, and $d\log (u):=du/u$. Moreover 
$H^1_{\fppf}(U,\mu_p)$ is identified with the set of 
regular differential forms on $U$ which are fixed by $C$ (cf. [Mi], III, 4).
\pop {.5}
\par
The Galois group $G_k$ acts in a natural way on 
$\mu _p$-torsors of $\overline U$. More precisely there exists a
canonical homomorphism:
$$(3')\ \ G_k\to \Aut H_{\fppf}^1(\overline U,\mu_p)$$
\par
Let $\sigma \in G_k$. Let $f:\overline V\to \overline U$ be a 
$\mu_p$-torsor, and let $f^{\sigma}:\overline V^{\sigma}\to 
\overline U$ be the $\mu_p$-torsor associated to $f$ via the above 
action $(3')$.  If $f$ is locally given by an equation $T_s^p=u_s$, 
then $f^{\sigma}$ is locally given by the equation $T_s^p=u_s^{\sigma}$,
where $u_s^{\sigma}$ is the image of $u_s$ via $\sigma$, 
which is an invertible function on $U_s^{\sigma}$.
In particular if $\omega:=du_s/u_s$ is the differential form associated to
$f$, then the differential form associated to $f^{\sigma}$ is
$\omega ^{\sigma}:=du_s^{\sigma}/u_s^{\sigma}$. The zeros of 
$\omega^{\sigma}$ contained in $\overline U$ are 
the images of those of $\omega$ via $\sigma$.
Moreover for a point $x_i\in \overline S$ we have 
$\cond _{x_i}(f)=\cond _{x_i^{\sigma}}(f^{\sigma})$ (resp. 
$\res _{x_i}(f)=\res _{x_i^{\sigma}}(f^{\sigma})$), where
$x_i^{\sigma}\in \overline S$ is the image of $x_i$ via $\sigma$.

\pop {.5}
\noindent
{\bf \gr 1.5. $\alpha _p$-Torsors.} 
\rm The following sequence is exact for the 
$\fppf$-topology:

$$(4)\ \  1\to \alpha _p\to \Bbb G_a@>{x^p}>>\Bbb G_a\to 1$$

Let $f:V\to U$ be an $\alpha _p$-torsor. There exists
an open covering $(U_s)_s$ of $U$ and regular fonctions $a_s\in \Gamma
(U_s,\Cal O_X)$, which are defined up to addition of a $p$-power, 
such that above $U_s$ the torsor $f$ is given by an equation 
$T_s^p=a_s$. The morphism $f$ is a finite radicial morphism 
of degree $p$, and in particular $V$ is homeomorphic to $U$. The 
differential form $\omega:=da_s$ is 
a global differential form on $X$, it is the 
{\bf differential form} associated to 
the torsor $f$. Let $\{y_j\}_{j\in Z}\subset U$ be the set of zeros of $w$ 
which are contained in $U$. Then as in 1.4 it is easy to see that
the only singularities of $V$ lie above the 
zeros of $w$ which are contained in $U$, and these singularities 
are unibranche.  For a point $x$ of $X$ we 
call the integer $\cond _{x}(f):=m_x:=-(\ord_x(\omega)+1)$ the 
{\bf conductor} of the torsor
$f$ at the point $x$. We define the {\bf residue} $h_x$ of the torsor $f$
at the point $x$ to be $0$. The differential form $\omega$
is annihilated by the Cartier operation $C$. In fact the following sequence
is exact on $U_{\et}$:
$$(5)\ 0\to \Cal O_U@>{F}>>
\Cal O_U@>{d}>>\Omega ^1_X@>{C}>>\Omega ^1_X\to 0$$
and $H^1_{\fppf}(U,\alpha_p)$ is identified with the set of 
regular differential forms on $U$ which are annihilated by $C$ 
(cf. [Mi], III, 4).

\pop {.5}
\par
The Galois group $G_k$ acts in a natural way on 
$\alpha_p$-torsors of $\overline U$. More precisely there exists a
canonical homomorphism:
$$(4')\ \ G_k\to \Aut H_{\fppf}^1(\overline U,\alpha_p)$$
\par
Let $\sigma \in G_k$. Let $f:\overline V\to \overline U$ be an 
$\alpha_p$-torsor, and let $f^{\sigma}:\overline V^{\sigma}\to 
\overline U$ be the $\alpha_p$-torsor associated to $f$ via the above 
action $(4')$. If $f$ is locally given by an equation $T_s^p=a_s$, 
then $f^{\sigma}$ is locally given by the equation $T_s^p=a_s^{\sigma}$,
where $a_s^{\sigma}$ is a regular function on $U_s^{\sigma}$.
In particular if $\omega:=da_s$ is the differential form associated to
$f$, then the differential form associated to $f^{\sigma}$ is
$\omega^{\sigma}:=da_s^{\sigma}$. The zeros of $\omega^{\sigma}$ contained in $\overline U$ are 
the images of those of $\omega$ via $\sigma$.
Moreover for a point $x_i\in \overline S$ we have 
$\cond _{x_i}(f)=\cond _{x_i^{\sigma}}(f^{\sigma})$ (resp. 
$\res _{x_i}(f)=\res _{x_i^{\sigma}}(f^{\sigma})$), where
$x_i^{\sigma}\in \overline S$ is the image of $x_i$ via $\sigma$.
\pop {1}

\epsfysize=5cm
\centerline{\epsfbox{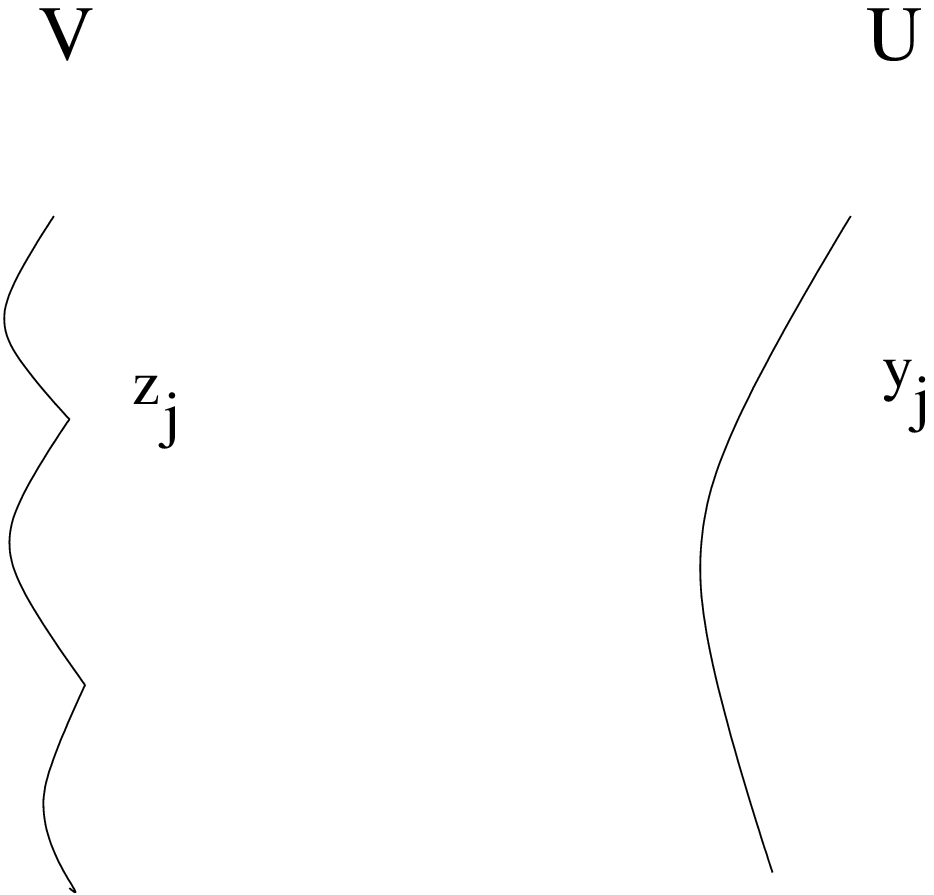}}

\pop {.5}
\par
\noindent
{\bf \gr 1.6. The group $H^1_{\fppf}(\ )_p$.}\rm
\par 
Let $Z$ be a connected scheme over a field $k$ of 
characteristic $p>0$.  Consider the following group:
$$(6)\ H^1_{\fppf}(Z)_p:=\oplus _{G_{k}}
H^1_{\fppf}(Z,G_{k})$$
where the sum is taken over all isomorphism classes of finite and flat 
commutative $k$-group schemes of rank $p$.
Let $\overline k$ be an algebraic closure of $k$, and let
$\overline Z:=Z\times _{\Spec k}\Spec {\overline k}$. Let $G_k$ be the 
Galois group of a sperable closure of $k$ contained in $\overline k$. We have:
$$H^1_{\fppf}(\overline Z)_p=H^1_{\fppf}(\overline Z,\mu_p)\oplus
H^1_{\et}(\overline Z,\Bbb Z/p\Bbb Z)\oplus 
H^1_{\fppf}(\overline Z,\alpha _p)$$ 
\par
The Galois group $G_k$ acts in a natural way on 
$H^1_{\fppf}(\overline Z)_p$ via its canonical action on 
$H^1_{\fppf}(\overline Z,\mu_p)$ 
(resp. on $H^1_{\fppf}(\overline Z,\alpha _p)$, and 
$H^1_{\fppf}(\overline Z,\Bbb Z/p\Bbb Z)$), and we have a canonical
homomorphism:
$$(5')\ G_k\to \Aut H^1_{\fppf}(\overline Z)_p$$ 

\par
In general for a scheme $Z$ with connected components
$\{{Z_i}\}_{i\in I}$, we define $H^1_{\fppf}(Z)_p:=\oplus _{i\in I}
H^1_{\fppf}(Z_i)_p$.

\pop {1}
\par
\noindent
{\bf \gr 1.7. The group $H^1_{\fppf}(\ )^{\kum}_p$.}
\rm
\par 
Let $X$ be a connected {\it semi-stable} curve over a field $k$ of 
characteristic $p>0$, with irreducible components
$\{X_i\}_{i\in I}$. Let $\{z_j\}_{j\in J}$ be the set of double points of 
$X$, and let $\{y_t\}_{t\in T}$ be a set of marked smooth points 
of $X$. 
Let $U:=X-\{\{z_j\}_{j\in J}\cup \{y_t\}_{t\in T}\}$. Then we define the group
$H^1_{\fppf}(U)^{\kum}_p$ as the subgroup of 
$H^1_{\fppf}(U)_p$ consisting of those elements $(f_s)_s$ with the following 
proprety, and which we call {\bf Kummerian mixed torsors} of $U$:
let $\{U_i\}_{i\in I}$ 
be the set of irreducible components of $U$. For each double point 
$z_j\in X_i\cap X_{i'}$ we impose that the following holds: there exists a 
component
$f_i:V_i\to U_i$ (resp. $f_{i'}:V_{i'}\to U_{i'}$) of $(f_s)_s$
which is a torsor under a finite and flat $k$-group scheme $G_{k,i}$ 
(resp. $G_{k,i'}$) of rank $p$, such that if $m_{i,j}$ (resp. 
$m_{i',j}$) is the conductor of the torsor $f_i$ (resp. $f_{i'}$) at 
the point $z_j$, and $h_{i,j}$ (resp. $h_{i',j}$) is the 
residue of the torsor $f_i$ (resp. $f_{i'}$) at the point $z_j$,
then $m_{i,j}+m_{i',j}=0$, and $h_{i,j}+h_{i',j}=0$.
It is not difficult to see that $H^1_{\fppf}(U)^{\kum}_p$ is indeed a 
subgroup of $H^1_{\fppf}(U)_p$ which is $G_k$-equivariant.

\pop {1}
\par
\noindent
{\bf \gr II. Degeneration of $\mu _p$-torsors from characteristic $0$ 
to characteristic $p>0$.}\rm \

\pop {.5}
\par
In this section we study the degeneration of $\mu_p$-torsors from 
zero to positive characteristics. We fix the following notations: 
$R$ is a complete discrete valuation ring of inequal 
characteristics, with residue characteristic $p>0$, and
which contains a primitive $p$-th root of unity $\zeta$. 
We denote by $K$ the fraction field of $R$, 
$\pi $ a uniformising parameter of $R$, $k$ the residue field of $R$, 
and $\lambda:=\zeta -1$. We also denote by $v_K$ the valuation of $K$ 
which is normalised by $v_K(\pi)=1$. Let $\overline K$ be an algebraic closure
of $K$, and let $G_K$ be the Galois group of $\overline K$ of $K$. 
Let  $\overline R$ be the integral closure of $R$ in $\overline K$ which
is a valuation ring, and let $\overline k$ be the residue field of 
$\overline R$ which is an algebraic closure of $k$. We denote by $G_k$ the 
Galois group of the separable 
closure of $k$ contained in $\overline k$. It is well known that there exists
a canonical exact sequence:
$$ (7)\ 0\to I_K\to G_K\to G_k\to 0$$
where $I_K$ is the inertia subgroup of $G_K$. Moreover the subgroup $I_K$ fits
in the following exact sequence:
$$0\to I_K^w\to I_K\to I_K^t\to 0$$
where $I_K^w$ (resp. $I_K^t$) is the wild part (resp. the tame part) 
of inertia. The subgroup $I_K^w$ of $I_K$ is a pro-$p$-group, 
and $I_K^t$ is canonically isomorphic to the prime to $p$-part 
$\Bbb Z(1)^{(p')}$ of the 
Tate twist $\Bbb Z(1)$ (cf. [Se]).

\pop {.5}
\par
\noindent
{\bf \gr 2.1. Torsors under finite and flat $R$-group schemes of rank $p$:
the group schemes $\Cal G^{n}$ and $\Cal H_n$.}
\rm\ (See also
[Oo-Se-Su] and [He]). 

\par
For any positif integer $n$ denote by 
$\Cal G^{n}_R:=\Spec R[x,1/(\pi ^nx+1)]$.
It is a commutative affine group scheme which has a generic fibre 
isomorphic to the multiplicative 
group $\Bbb G_m$ and whose special fibre is the additive group $\Bbb G_a$
(cf. [Oo-Se-Su] for more details).
For $0< n\le v_K(\lambda)$ the polynomial $((\pi ^nx+1)^p-1)/\pi ^{pn}$ has 
coefficients in $R$ and defines a group scheme homomorphism $\varphi _n:
\Cal G^n_R\to \Cal G^{np}_R$ given by:
$$R[x,1/(\pi ^{pn}x+1)]\to R[x,1/(\pi ^nx+1)]$$
$$x\to ((\pi ^nx+1)^p-1)/\pi ^{pn}$$
\par
The isogeny $\varphi _n$ has degree $p$, its kernel $\Cal H_{n,R}$
is a finite and flat $R$-group scheme of rank $p$, its generic fibre is 
isomorphic 
to $\mu _p$, its special fibre is either the radiciel group $\alpha _p$
if $0<n< v_K(\lambda)$, or the \'etale group $\Bbb Z/p\Bbb Z$ if 
$n=v_K(\lambda)$. 
\par
Let $\alpha ^{(n)}:\Cal G^n_R\to \Bbb G_{m,R}$ be the 
group schemes homomorphism given by:
$$R[u,u^{-1}]\to R[x,1/(\pi ^nx+1)]$$
$$u\to \pi ^nx+1$$
\par
The following is a commutative diagramm of exact sequences 
for the $\fppf$-topology:

$$
\CD
0@>>>    \Cal H_{n,R}    @>>>     \Cal G^n_R   @>{\varphi _n}>>\Cal G^{np}_R@>>> 0\\
 @.        @VVV            @V\alpha^{(n)}VV              @V\alpha ^{(np)}VV \\
0      @>>>    \mu _{p,R}   @>>>   \Bbb G_{m,R}@> x^p>> \Bbb G_{m,R}@>>> 0
\endCD
$$

\par
The upper exact sequence above has a Kummer exact sequence (1) as 
generic fibre, and the exact sequence (4) if $0<n< v_K(\lambda)$, or an 
exact sequence of Artin-Schreir type (2) if $n=v_K(\lambda)$, as 
a special fibre. 
\par
Let $\Cal U$ be an $R$-scheme, and let $f:\Cal V\to \Cal U$ be 
a torsor under $\Cal H_{n,R}$. Then there exists an open covering
$(\Cal U_i)$ of $\Cal U$ and regular fonctions $u_i\in \Gamma 
(\Cal U_i,\Cal O_{\Cal U})$, such that $\pi ^{np}u_i+1$ is defined up to 
multiplication by a $p$-power, and the torsor $f$ is defined above $\Cal U_i$
by an equation $T_i^p=(\pi ^{n}T'_i+1)^p=\pi ^{np}u_i+1$, where $T_i$ 
and $T'_i$ are indeterminates.

\pop {.5}
\par
\noindent
{\bf \gr 2.2. The Galois action.}\rm 
\par
Let $\Cal G'$ be a finite 
and flat commutative $R$-group scheme of rank 
$p$, and let ${\Cal G}:=\Cal G'\times _R{\overline R}$ which is a 
finite and flat $\overline R$-group scheme of rank $p$. We denote by 
${\Cal G}_{\overline k}:=
{\Cal G}\times _{\overline R}{\overline k}$ the special fibre of 
${\Cal G}$ which is either $\Bbb Z/p\Bbb Z$, $\mu_p$, or 
$\alpha _p$. Let $\Cal U$ be an $R$-scheme, and let 
$\overline {\Cal U}:=\Cal U\times _R{\overline R}$. 
The Galois group $G_K$ acts in a canonical way on
torsors, for the $\fppf$-topology, above $\overline {\Cal U}$ 
under the group scheme ${\Cal G}$. More precisely there exists a canonical 
homomorphism:
$$(6')\ \ G_K\to \Aut H_{\fppf}^1(\overline {\Cal U},{\Cal G})$$
\par
The above action $(6')$ is compatible with the action of $G_k$
on torsors above $\overline {\Cal U}_{\overline k}:=
\overline {\Cal U}\times _{\overline R}{\overline k}$ under the group scheme
${\Cal G}_{\overline k}$. More precisely we have the following 
canonical commutative diagramm:

$$
\CD
G_K  @>>>   \Aut H_{\fppf}^1(\overline {\Cal U},{\Cal G})  \\
     @VVV            @VVV \\
G_k    @>>> \Aut H_{\fppf}^1(\overline {\Cal U}_{\overline k},{\Cal G}_{\overline k})
\endCD
$$

\pop {.5}
\par
\noindent
{\bf \gr 2.3. Degeneration of $\mu _p$-torsors.}\rm \ 
\par
In what follows let $X$ be a formal 
$R$-scheme of finite type which is normal and flat over $R$.
We assume that $X$ is {\bf smooth of relative dimension 1}. Let 
$X_K:=X\times _R K$ (resp. $X_k:=X\times_R k$) be the generic (resp. 
special) fibre of $X$. By the generic fibre of $X$ we mean the associated 
$K$-rigid space. Let $\eta$ be the generic point of the special 
fibre $X_k$, and let $\Cal O_{\eta}$  be the local ring of $X$ at $\eta $, 
which is a discrete valuation ring with fractions fields $K(X)$ the 
function field of $X$. Let $f_K:Y_K\to X_K$ be a  
$\mu _p$-torsor, and let 
$K(X)\to L$ be the corresponding extension of function fields. 
One has the following:

\pop {.5}
\par
\noindent
{\bf \gr 2.4. Proposition.}\rm\ {\sl Assume that the ramification index
above $\Cal O_{\eta}$ in the extension $K(X)\to L$ equals $1$. Then
the torsor $f_K:Y_K\to X_K$ extends 
to a torsor $f:Y\to X$ under a finite and flat $R$-group scheme of rank $p$,
with $Y$ normal. Let 
$\delta$ be the 
degree of the different above $\eta$ in the extension $K(X)\to L$. Then the 
following cases occur:    
\par
a ) $\delta = 0$ in which case $f$ is a torsor under the group scheme 
$\Cal H_{v_K(\lambda),R}$, and $f_k$ is then an \'etale torsor under 
$\Bbb Z/p\Bbb Z$.
\par
b ) $0<\delta<v_K(\lambda)$ in which case $\delta = v_K(\lambda)-n(p-1)$
for a certain integer $n\ge 1$,  
$f$ is a torsor under $\Cal H_{n,R}$ and $f_k$ is a radicial torsor under $\alpha _p$.
\par
c ) $\delta =v_K(\lambda)$, $f$ is a torsor under $\mu _p$ and $f_k$ is also a 
radicial torsor under $\mu _p$.}

\pop {.5}
\par
Note that starting from a torsor $f_K:Y_K\to X_K$ as in 2.3, the condition 
that the ramification index above $\Cal O_{\eta}$ equals 1 is alaways
satisfied after eventually a finite ramified extension of $R$ (cf. [Ep]).

\pop {.5}
\par
\noindent
{\bf \gr Proof.}\rm\ The above result is classically known in the case where 
$X$ is the formal 
spectrum of a complete discrete valuation ring (cf. [Hy], 
note that this is the only case we will need in 2.4). 
It has been treated in the
case where $X$ is factoriel in [Ma-Gr] and [He]. 
We explain the proof in the general case. 
By hypothesis there exists an open covering $(U_i)_i$ of $X$ and units
$u'_i\in \Gamma (U_{i,K},\Cal O_X)^*$ which are defined up to multiplication 
by $p$-powers, where $U_{i,K}:=U_i\times _R K$, and such that the torsor $f_K$
is defined above $U_{i,K}$ by an equation $T_i^p=u'_i$. Moreover $u'_i=
\pi ^{p\alpha_i}u_i$ where $u_i\in \Gamma (U_{i},\Cal O_X)^*$ is a unit 
(the power of $\pi$ is a $p$-multiple since the ramification index above 
$\eta$ equals $1$, and the fact that $u_i$ is a unit can be cheked easily 
locally and uses the fact that $X$ is a smooth curve). If for a 
given $i$ the image $\bar u_i$
of $u_i$ in $\Cal O_X(U_i)/\pi\Cal O_X(U_i)$ is not a $p$-power then this 
is the case for each $i$ since $u_i=a_{i,j}^pu_j$ in $\Cal O_X(U_i\cap U_j)$. 
And in this case the class of the cocycle $(a_{i,j})_{i,j}$ defines a non 
trivial $\mu_p$-torsor $f:Y\to X$ of $X$ which extends $f_K$ and induces a 
$\mu _p$-torsor of $X_k$. Otherwise $\bar u_i$ is a $p$-power for
each $i$. Then as in [Gr-Ma], [He], and after multiplying $u_i$ by a suitable
$p$-power, which does not change the $\mu _p$-torsor $f_K$, one can assume 
that $f_i=1+\pi ^{pt_i}a_i$, and the image $\bar a_i$ of $a_i$ in
$\Cal O_X(U_i)/\pi\Cal O_X(U_i)$ is not a $p$-power, in which case 
$\delta _i:=v(\lambda)-t_i(p-1)$ is the degree of the different 
above the generic point 
$\eta$ of $U_{i,k}:=U_i\times _R k$. In particular $n:=t_i$ is constant for 
all $i$. The $\Cal H_{n,R}$-torsor $f:Y\to X$ above $X$, which is given 
above $U_i$ by the 
equation $T_i^p=1+\pi ^{pn}a_i$, extends then the $\mu _p$-torsor $f_K$.

\pop {.5}
\par
\noindent
{\bf \gr 2.5.}\rm \ With the same notations as in 2.3 let
$X_{\overline K}:=X\times _R{\overline K}$,  
$X_{\overline R}:=X\times _R{\overline R}$, and let $X_{\overline k}:=
X\times _R{\overline k}$. Hier we consider $X_{\overline K}$ as a rigid
space. Then it follows easily from the proof of 
2.4, that there exists canonical {\bf specialisation} group homomorphisms:
$$H^1_{\et}(X_{\overline K},\mu_p)\to H^1_{\fppf}(X_{\overline R},\mu_p)
\to H^1_{\fppf}(X_{\overline k},\mu_p)$$ 
whose composite gives the following:

\pop {.5}
\par
\noindent
{\bf \gr 2.5.1. Proposition.}\rm\ {\sl There exists a canonical 
specialisation group homomorphism:
$$\Sp: H^1_{\et}(X_{\overline K},\mu_p)\to 
H^1_{\fppf}(X_{\overline k},\mu_p)$$ 
which is $G_K$-equivariant.}
\pop {.5}
\par
Note that the $\mu_p$-torsors of $X_{\overline K}$ which 
belongs to $\Ker (\Sp)$ correspond to thoses torsors which induce in reduction
$\alpha _p$ or $\Bbb Z/p\Bbb Z$-torsors. For the proof of the second 
assertion of 2.5.1 cf. 2.7. Also hier we mean the action of $G_K$ on  
$H^1_{\fppf}(X_{\overline k},\mu_p)$ via its canonical quotient:
$G_K\to G_k$. 
\pop {.5}
\par
\noindent
{\bf \gr 2.6. Compatibility of degeneration with the Galois 
action.}
\par
\rm\ Let $X$ be a formal $R$-scheme of finite type which is normal and flat 
over $R$. We suppose that $X$ is {\bf smooth of relative dimension one}. Let 
$X_K:=X\times _R K$ (resp. $X_k:=X\times_R k$) be the generic 
(resp. special) fibre of $X$. By the generic fibre of $X$ we mean the 
associated $K$-rigid space. We have the Galois action of $G_K$ 
on $\mu_p$-torsors of $X_{\overline K}:=X\times _K \overline K$ (cf. $(1')$):
$$\ G_K\to H^1_{\et}(X_{\overline K},\mu_p)$$
\par
As a direct consequence of the result in 2.4 one sees that there exists a 
canonical {\bf specialisation} set theoretical map 
(this is not a group homomorphism):
$$\ \Sp:H^1_{\et}(X_{\overline K},\mu_p)\to H^1_{\fppf}(X_{\overline k})_p$$
Where $X_{\overline k}:=X\times _{\Spec k}\Spec {\overline k}$. Moreover the 
action of $G_k$ on $H^1_{\fppf}(\overline X)_p$ induces  
an action of $G_K$ on $H^1_{\fppf}(\overline X)_p$ via the canonical surjectif 
homomorphism: $G_K\to G_k$. We have the following:

\pop {.5}
\par
\noindent
{\bf \gr 2.7. Theorem.}\rm \ {\sl There exists a canonical 
specialisation map 
$\Sp:H^1_{\et}(X_{\overline K},\mu_p)\to H^1_{\fppf}(X_{\overline k})_p$ 
which is $G_K$-equivariant.}

\pop {.5}
\par
\noindent
{\bf \gr Proof.}\rm \ We explain the argument of the proof. We may assume that 
$\Pic (X_{\overline K})=0$. Let $f_{\overline K}:Y_{\overline K}\to 
X_{\overline K}$
be a non-trivial $\mu_p$-torsor. Then $f_{\overline K}$ is given by an 
equation $T^p=u_{\overline K}$, where $u_{\overline K}$ 
is an invertible function on $X_{\overline K}$. By 2.4 this torsor 
extends to a torsor $f_{\overline R}:Y_{\overline R}\to X_{\overline R}$ 
under a finite
and flat ${\overline R}$-group scheme of rank $p$, where 
$X_{\overline R}:=X\times _R {\overline R}$ and ${\overline R}$ is the integral
closure of $R$ in $\overline k$. This torsor induces on the special fibres
a torsor $f_{\overline k}:Y_{\overline k}\to X_{\overline k}$ under a finite
and flat ${\overline k}$-group scheme of rank $p$, where 
$X_{\overline k}:=X\times _k {\overline k}$, and by definition 
$f_{\overline k}$ is the non zero component of the image of 
$f_{\overline K}$ in 
$H^1_{\fppf}(X_{\overline k})_p$ via the map $\Sp$. 
We treat the case where $f_{\overline R}$ is an \'etale torsor. The other 
cases are treated similarly. Then it follows from 2.4 that one can choose
the function  $u_{\overline K}$ of the form $1+\lambda ^p
a$, where $a$ is a regular function on 
$X_{\overline R}$. Moreover the torsor $f_k$ is given by the equation
$t^p-t=\overline a$, where $\overline a$ is the image of
$a$ in $\Gamma (X_{\overline k},\Cal O_{X_{\overline k}})$.
 Let $\sigma \in G$. Then the torsor ${f_{\overline R}}^{\sigma}$, image
of $f_{\overline R}$ under $\sigma$, is given 
by the equation $1+\lambda^p{a}^{\sigma}$,
where ${a}^{\sigma}$ is the transform of 
${a}$ under $\sigma$,
its special fibre is the torsor given by the equation 
$t^p-t={\overline a}^{\overline {\sigma}}$, where 
${\overline a}^{\overline {\sigma}}$ is the transform  of 
${\overline a}$ under ${\overline {\sigma}}$ which is the image of 
$\sigma$ in $G_k$. But then this torsor is nothing else but
the transform ${f_{\overline k}}^{\overline {\sigma}}$
of $f_{\overline k}$ under ${\overline {\sigma}}$.

\pop {.5}
\par
\noindent
{\bf \gr III. Degeneration of $\mu _p$-torsors on the boudaries of 
formal fibres.}

\pop {.5}
\par
In this section we explain the degeneration of $\mu_p$-torsors on the 
boundary 
\newline
$\Spf R[[T]]\{T^{-1}\}$ of formal fibres of formal $R$-curves, where
$R[[T]]\{T^{-1}\}$ is the ring of formal power
series $\sum_{i\in \Bbb Z} a_iT^{i}$ with $\lim _{i\to -\infty} \vert a_i
\vert =0$, and where $\vert\ \vert$ is an absolute value of $K$ associated 
to its valuation. The following result will be used in [Sa] in order 
to prove a formula for vanishing cycles.

\pop {.5}
\par
\noindent
{\bf \gr 3.1. Proposition.}\rm \ {\sl Let $A:=R[[T]]\{T^{-1}\}$ 
(cf. the definition above), 
and let 
$f:\Spf B\to \Spf A$ be a non trivial torsor under a finite flat 
$R$-group scheme $G$ of 
rank $p$. let $\delta$ be the degree of the different in the above 
extension. We assume that the ramification index of this extension
equals $1$. Then the following cases occur:
\par
a ) $\delta =0$, the torsor $f$ is \'etale and is given, after eventually
a finite extention of $R$,  by an equation 
$Z^p={\lambda ^p}T^{m}+1$, where $m$ is a negative integer prime to $p$, for 
a suitable choice of the parameter $T$ of $A$.
\par
b ) $\delta=v_K(\lambda)$, $f$ is a torsor under $\mu _p$, and two cases 
can occur:
\par
b-1 ) For a suitable choice of the parameter $T$ of $A$, and 
after eventually a finite extention of $R$, the torsor $f$ is given, 
by an equation $Z^p=T^{h}$ where $h\in \Bbb F_p^*$.  
\par
b-2 ) For a suitable choice of the parameter $T$ of $A$, and 
after eventually a finite extention of $R$, the
torsor $f$ is given by an equation $Z^p=1+T^m$ where $m$ is a positive 
integer prime to $p$.
\par
c ) $0 <\delta <v_K(p)$ in which case $f$ is a torsor under $\Cal H_{n,R}$,
and $\delta = v_K(p
)-n(p-1)$. For a suitable choice of the 
parameter $T$, and 
after eventually a finite extention of $R$, the torsor $f$ is given
by an equation $Z^p=1+\pi ^{np}T^m$, with $m\in \Bbb Z$ is prime to $p$}.

\pop {.5}
\par
The proof of the above lemma is easy. It is a direct 
consequence of 2.4 which in this case is due 
to Hyodo (cf. [Hy] 2.16). This lemme is also stated in [He] in a slight 
different form in terms of automorphisms of boundaries of formal fibres.

\pop {.5}
\par
\noindent
{\bf \gr 3.2. Definition.}\rm \ With the same notations as in 3.1,
we define the {\bf type of reduction} (or the degeneration type) 
of the torsor $f$ to be $(G_k,m,h)$,
where $G_k:=G\times _R k$ is the special fibre of $G$, $m$ is 
the conductor associated to the torsor $f_k:\Spec B/\pi B\to \Spec A/\pi A$, 
and $h\in \Bbb F_p$ its residue (cf. 1.3, 1.4 and 1.5). With the same 
notations as in 2.4 the degeneration type is $((\Bbb Z/p\Bbb Z),m,0)$ in
case a, $(\mu _p,0,h,)$ in case b-1, $(\mu _p,-m,0)$ in case b-2 and
$(\alpha _p,-m,0)$ in case c.

\pop {.5}
\par
\noindent
{\bf \gr 3.3. Remark.}\rm\ The above corollary implies in particular 
that given two torsors above $A:=R[[T]]\{T^{-1}\}$, under a finite 
flat $R$-group scheme of rank $p$, and which have the same type of 
reduction as defined in 3.2, and which have the same 
degree of the different, then after eventually ``adjusting'' 
the Galois action on the Kummer generators of these two covers, and
after eventually a finite extension of $R$, one can find 
a Galois-equivariant isomorphism between them. Also note that the above 
lemma can be easily adapted to the rigid setting in order to describe 
torsors under group schemes of rank $p$ above the boundaries of formal 
fibres of curves in rigid geometry as well as their degeneration type.

\pop {.5}
\par
\noindent
{\bf \gr IV. Lifting of rank $p$ torsors from positive 
to zero characteristic.}\rm\ 
\par
In this section we use the same notations as in 2.3.  
Let $X$ be a formal $R$-scheme of finite type, which is normal and flat 
over $R$. We assume that $X$ is smooth of relative dimension one. 
Let $X_k:=X\times_R k$ be the special fibre of $X$. 
Also let $X_{\overline R}:=X\times _R\overline R$, 
and let $X_{\overline k}:=X\times _k{\overline k}$.

\pop {.5}
\par
\noindent
{\bf \gr 4.1. Lifting of rank $p$ \'etale torsors.}
\rm

\par
Let $f_k:Y_k\to X_k$ be an \'etale $\Bbb Z/p\Bbb Z$-torsor. Then 
by the theorems of lifting of \'etale covers (cf. [Gr]) $f_k$ can 
be lifted, uniquely up to isomorphism, to an \'etale torsor 
of $X$. More precisely if the torsor $f_k$ is locally given by an equation
$t^p-t=\bar a$, where $\bar a$ is a regular function, then the equation
$((\lambda T+1)^p-1)/\lambda ^p=a$ where $a$ is a regular function which equals
$\bar a$ mod $\pi$ gives locally a lifting of the torsor $f_k$. In particular
if $R^{\nr}$ denotes the maximal unramified extension of $R$ contained in
$\overline R$, then all \'etale torsors under a finite and flat 
$\overline R$-group scheme of rank $p$ are already defined over $R^{\nr}$, and 
the Galois action of $G_K$ on those torsors factors through 
the Galois group of $K^{\nr}$: the fractions field of $R^{\nr}$, 
over $K$, which is canonically isomorphic to $G_k$. More precisely we have 
$G_k$-equivariant canonical isomorphisms:
$$H^1_{\et}(X_{\overline k},\Bbb Z/p\Bbb Z)\simeq 
H^1_{\fppf}(X_{\overline R},\Cal H_{\lambda,\overline R})\simeq
H^1_{\fppf}(X_{R^{\nr}},\Cal H_{\lambda,R^{\nr}})$$
where $X_{R^{\nr}}:=X\times _R R^{\nr}$,  
$\Cal H_{\lambda,R^{\nr}}:=\Cal H_{\lambda,R}\times _R{R^{\nr}}$,
and $\Cal H_{\lambda,\overline R}:=\Cal H_{\lambda,R}\times _R{\overline R}$
\par
The theorem of lifting of \'etale covers above can also be interpreted as
follows:

\pop {.5}
\par
\noindent
{\bf \gr 4.1.1. Lemma.}\rm\ {\sl There exists
a canonical $G_K$-equivariant injective group homomorphism:
$$H^1_{\et}(X_{\overline k},\Bbb Z/p\Bbb Z)\to H^1_{\et}
(X_{\overline K},\mu_p)$$ 
\par
Moreover the Galois group $G_K$ acts on the subgroup
$H^1_{\et}(X_{\overline k},\Bbb Z/p\Bbb Z)$ via its canonical
quotient $G_{{K^{\nr}}/K}$: the Galois group of $K^{\nr}$ over $K$,
which is canonically isomorphic to $G_k$.} 
\par
In other words the inertia 
subgroup $I_K$ of $G_K$ acts trivially on the subgroup 
\newline
$H^1_{\et}(X_{\overline k},\Bbb Z/p\Bbb Z)$.
\pop {.5}
\par
\noindent
{\bf \gr 4.2. Lifting of rank $p$ radiciel torsors.}
\rm\ 
\pop {.5}
\par
\noindent
{\bf \gr 4.3. Proposition.}\rm \ {\sl With the same notations as in IV above, 
let $f_k:Y_k\to X_k$ be a torsor above $X_k$ under a 
finite and flat $k$-group scheme $G_k$ of rank $p$ which is radicial. 
Assume that $X:=\Spf A$ is {\bf affine} and that
$H^1_{\et}(X_k,\Bbb G_m)=
H^1_{\fppf}(X_k,\Bbb G_m)=\Pic (X_k)=0$ is {\bf trivial}. 
Then the torsor $f_k$ can be lifted, non uniquely, 
after eventually a finite ramified extention of $R$, to a torsor 
$f:Y\to X$ under a finite and flat $R$-group scheme of rank $p$.} 

\pop {.5}
\par
More precisely we have the following two cases:
\pop {.5}
\par
\noindent
{\bf case 1: $G_k=\mu_ p$.}\ In this case the torsor $f_k$ is given by an 
equation $t^p=\bar u$, where $\bar u\in A_k:=A/\pi A$ is a unit which is 
uniquely determined up to multiplication by a $p$-power. 
Let $u\in A$ be an invertible function such that $u$ equals 
$\bar u$ mod $\pi$. Then the equation $T^p=u$ defines a $\mu_p$-torsor
$f:Y\to X$ above $X$ which lifts the $\mu_p$-torsor $f_k$. The above lifting 
is however not unique. More precisely for $n\le v_K(\lambda)$ a 
positif integer,
the equation ${T'}^p=u(1+\pi^{pn}u')$, where $u'\in  A$ defines 
another lifting $f':Y'\to X$ of $f_k$. Moreover $f$ and $f'$ are not 
isomorphic $\mu_p$-torsors since $(1+\pi^{pn}u')$ is not a $p$-power in $A$.
All possible liftings of the torsor $f_k$ are in fact defined up to elements
in $1+\pi ^m A$, for $n\le pv_K(\lambda)$ a positif integer. Moreover
in this case the composite of the canonical homomorphisms: 
$H^1_{\fppf}(X_{\overline K},\mu_p)\to H^1_{\fppf}(X_{\overline R},\mu_p)
\to H^1_{\fppf}(X_{\overline k},\mu_p)$ is surjective, and we have a 
canonical commutative diagramm:
$$
\CD
G_K  @>>>   \Aut H_{\fppf}^1(X_{\overline K},\mu_p)  \\
   @VVV            @VVV \\
G_K  @>>>   \Aut H_{\fppf}^1(X_{\overline R},\mu_p)  \\
     @VVV            @VVV \\
G_k    @>>> \Aut H_{\fppf}^1(X_{\overline k},\mu _p)
\endCD
$$

\par
\noindent
{\bf case 2: $G_k=\alpha _p$.}\ In this case the torsor $f_k$ is given by an 
equation $t^p=\bar a$, where $\bar a\in A_k:=A/\pi A$ is a function which is 
uniquely determined up to addition of a $p$-power. Let $n< v_K(\lambda)$ be
a positif integer. Such an integer exists if and only if $v_K(\lambda)>1$.
Note that for any positif integer $n$, one can perform a totally wildely 
ramified extention of $R$ such that the condition $n< v_K(\lambda)$ 
is satisfied. Assume that such an $n$ as above exists. Let $a\in A$ 
be a function which equals $\bar a$ mod $\pi$. Then the equation 
$T^p=1+\pi^{np}a$ defines an $\Cal H_{n,R}$-torsor $f:Y\to X$ above 
$X$ which lifts the $\mu_p$-torsor $f_k$. The above lifting 
is however not unique. More precisely for $m$ a positif
integer such that $m\le v_K(\lambda)$, and $pn<pm$,
the equation ${T'}^p=(1+\pi^{pn}a)(1+\pi^{pm}a')$, where $a'\in  A$ defines 
another lifting $f':Y'\to X$ of $f_k$. Moreover $f$ and $f'$ are not 
isomorphic $\Cal H_{n,R}$-torsors since $(1+\pi^{pm}a')$ is not a $p$-power 
in $A$. 
\par
Suppose there exists a positif integer $n < v_K(\lambda)$, and let
$\Cal H_{n,\overline R}:=\Cal H_{n,R}\times _R{\overline R}$ which
is a finite and flat commutatif $\overline R$-group scheme of rank $p$, 
whose special fibre is the radicial $\overline k$-group scheme 
$\alpha _p$. Then the canonical 
homomorphism: $H^1_{\fppf}(X_{\overline R},\Cal H_{n,\overline R})\to
H^1_{\fppf}(X_{\overline k},\alpha_p)$ is surjective, and we have a 
commutative diagramm:
$$
\CD
G_K  @>>>   \Aut H_{\fppf}^1(X_{\overline R},\Cal H_{n,\overline R})  \\
     @VVV            @VVV \\
G_k    @>>> \Aut H_{\fppf}^1(X_{\overline k},\alpha _p)
\endCD
$$

\pop {.5}
\par
\noindent
{\bf \gr 4.4.}\rm \ Let $G_k$ be a radicial finite and flat $k$-group 
scheme of rank $p$. Let $f_k:Y_k\to X_k$ be a torsor above $X_k$ under
$G_k$. We say that this torsor is {\bf admissible}
if, eventually after a finite extension of $R$, there exists a finite and 
flat commutatif $R$-group scheme $G$, whose special fibre is isomorphic 
to $G_k$, and a torsor $f:Y\to X$ above $X$ under $G$ which lifts $f_k$.
By the considerations above we have seen that if $X$ is affine and 
$\Pic (X)=0$, then every torsor $f_k:Y_k\to X_k$ above $X_k$ under 
$G_k$ is admissible. This is for example the case if $X$ is a formal
affine open of the formal $R$-projective line. Another important example 
of admissible torsors is the case where $X$ is a proper and smooth formal 
$R$-curve. More precisely we have the following: 
\pop {.5}
\par
\noindent
{\bf \gr 4.5. Proposition.}\rm\ {\sl Let $X$ be a {\bf proper} and 
{\bf smooth} formal $R$-curve. Let $G_k$ be  a commutative finite and flat 
$k$-group scheme of rank $p$ which is radicial. Then every torsor 
$f_k:Y_k\to X_k$ above $X_k$ under the group scheme $G_k$ is admissible.}
\pop {.5}
\par
\noindent
{\bf \gr Proof.}\rm\ Let $J:=\Pic ^0(X)$ be the jacobian of $X$
which is an abelian formal scheme. Let $G$ be a finite and flat commutatif
$R$-group scheme, and let $G'$ be its Cartier dual. 
Then by [Ra] there exists a canonical isomorphism
$H^1_{\fppf}(X,G)\simeq \Hom (G',J)$. Let $G_k'$ be the Cartier dual
of $G_k$ which by the geometric abelian class field theory is a subgroup  
of $J_k:=J\times _Rk$ the special fibre of $J$. To lift then the torsor
$f_k$ is equivalent by the above isomorphism to lift the subgroup
$G_k'$ of $J_k$ to a subgroup $G'$ of $J$ which is finite and flat over 
$R$ of rank $p$, and this is indeed alaways possible.

\pop {.5}
\par
\noindent
{\bf \gr 4.6. Remark.}\rm\ I do not know whether for any 
formal affine $R$-curve $X$, and a finite flat $k$-group scheme $G_k$ 
of rank $p$, any $G_k$-torsor of $X_k$ is admissible?.

\pop {.5}
\par
\noindent
{\bf \gr V. A canonical filtration on $H^1_{\et}(\ ,\mu_p)$.}
\rm
\par
In this section we use the same notations as in III. Let $X$ be an 
$R$-scheme of finite type, which is {\bf smooth}, {\bf proper} and flat 
over $R$ of {\bf relative dimension one}. Let 
$X_k:=X\times_R k$ be the special fibre of $X$. Also let 
$X_{\overline K}:=X\times _R\overline K$, and let $X_{\overline k}:=
X\times _k{\overline k}$. We will establish
in this section the existence of a canonical $G_K$-equivariant 
{\bf filtration} $\Fil:=(\Fil _i)_i$ on the group 
$H^1_{\et}(X_{\overline K},\mu_p)$.

\par
Recall first (cf. 2.5) that we have a canonical specialisation group
homomorphism which is $G_K$-equivariant:
$$\Sp: H^1_{\et}(X_{\overline K},\mu_p)\to 
H^1_{\fppf}(X_{\overline k},\mu_p)$$ 
and the kernel of $\Sp$ consist of those $\mu_p$-torsors of $X_{\overline K}$
which induces in reduction an $\alpha_p$ or an \'etale torsor of 
$X_{\overline k}$.
\par
The group $H^1_{\et}(X_{\overline K},\mu_p)$ is finite 
since $X$ is proper. In particular all $\mu_p$-torsors of $X_{\overline K}$
are defined over a finite extension of $K$. More precisely there exists a 
finite extension $K'$ of $K$ with the following properties:

\par
(a)\ All $\mu_p$-torsors of $X_{\overline K}$ are defined over $K'$.
\par
(b)\ Let $f_{K'}:Y_{K'}\to X_{K'}:=X\times _R K'$ be a $\mu_p$-torsor, 
and let $Y'$ be the integral closure of $X'$ in $Y_{K'}$, where
$X':=X\times _R R'$, and $R'$ is the integral closure of $R$ in $K'$. 
Then the special fibre of $Y'$ is reduced.

\par
Moreover there exists such a finite extension $K'$ of $K$ which 
is minimal for the above properties. Let $\pi'$ be a uniformiser of $K'$,
and let $v_{K'}$ be the valuation of $K'$ which is normalised by 
$v_{K'}(\pi ')=1$. It follows from 2.3.4 that if  
$f_{K'}:Y_{K'}\to X_{K'}:=X\times _R K'$ is a non trivial $\mu_p$-torsor, 
then it extends to a torsor $f:Y'\to X':=X\times _R R'$ under a finite 
commutative and flat $R'$-group scheme. Let $\delta _f:=\delta _{{f_K}}$
be the degree of the different, in the above torsor $f$, above the 
generic point of the special fibre of $X'$. We have two cases: either
$f$ is a non trivial $\mu_p$-torsor, or a torsor under the
group scheme $\Cal H_{n,R'}$ for $0<n\le v_{K'}(\lambda)$, the case
$n=v_{K'}(\lambda)$ corresponds to the case of an \'etale torsor.
In the first case we have $\delta _f=v_{K'}(p)$, and in the second case we 
have $\delta _f=v_{K'}(p)-n(p-1)$.

\pop {.5}
\par
\noindent
{\bf \gr 5.1. Definition.}\rm\ With the same notations as above we define 
the following decreasing filtration  $\Fil:=(\Fil _n)_{0\le n\le v_{K'}(p)}$
of the group $H^1_{\et}(X_{\overline K},\mu_p)$: 
$\Fil _n H^1_{\et}(X_{\overline K},\mu_p)$ is the subgroup of
$H^1_{\et}(X_{\overline K},\mu_p)$ consisting of those, isomorphism classes,
of $\mu_p$-torsor $f_{K'}:Y_{K'}\to X_{K'}:=X\times _R K'$ with
$\delta _{f_{K'}}\le v_{K'}(p)-n(p-1)$. In particular we have: 
$\Fil _0 H^1_{\et}(X_{\overline K},\mu_p)=
H^1_{\et}(X_{\overline K},\mu_p)$, and 
$\Fil _{v_{K'}(\lambda)} H^1_{\et}(X_{\overline K},\mu_p)\simeq
H^1_{\et}(X_{\overline k},\Bbb Z/p\Bbb Z)$.

\pop {.5}
\par
The following proposition is easily verified:

\pop {.5}
\par
\noindent
{\bf \gr 5.2. Proposition.}\rm\ {\sl The filtration 
$\Fil:=(\Fil _n)_{0\le n\le v_{K'}(p)}$ is a decreasing filtration 
by subgroups of $H^1_{\et}(X_{\overline K},\mu_p)$. Moreover
this filtration is $G_{K}$-equivariant.}

\pop {.5}
\par
\noindent
{\bf \gr 5.3.}\rm \ With the same notations as above, let
$g$ be the genus of $X_{\overline K}:=X\times _R{\overline K}$, and let 
$r$ be the $p$-rank
of $X_{\overline k}:=X\times _R {\overline k}$ which equals 
$\dim _{\Bbb Z/p\Bbb Z}H^1_{\et}
(X_{\overline k},\Bbb Z/p\Bbb Z)$, and which is an integer smaller or
equal to $g$. Then it is well known that:
$H^1_{\et}(X_{\overline K},\mu _p)\simeq 
(\Bbb Z/p\Bbb Z)^{2g}$, 
$H^1_{\fppf}(X_{\overline k},\alpha _p)\simeq 
(\Bbb Z/p\Bbb Z)^{2(g-r)}$, and  $H^1_{\fppf}(X_{\overline k},\mu _p)
\simeq (\Bbb Z/p\Bbb Z)^r$. Moreover in this case the canonical 
$G_K$-equivariant homomorphism:
$$\Sp: H^1_{\et}(X_{\overline K},\mu_p)\to 
H^1_{\fppf}(X_{\overline k},\mu_p)$$ 
is surjectif. Moreover we have: 
$\Fil _0 H^1_{\et}(X_{\overline K},\mu_p)=H^1_{\et}(X_{\overline K},\mu_p)$,
and 
$\Fil _1 H^1_{\et}(X_{\overline K},\mu_p)=\Ker (\Sp)$. In particular
$\Fil _1 H^1_{\et}(X_{\overline K},\mu_p)\simeq (\Bbb Z/p\Bbb Z)^{2g-r}$.
Also: $\Fil _{v_{K'}(\lambda)} H^1_{\et}(X_{\overline K},\mu_p)\simeq
H^1_{\et}(X_{\overline k},\Bbb Z/p\Bbb Z)
\simeq (\Bbb Z/p\Bbb Z)^r$.
\pop {.5}
\par
In the case where $X_{\overline k}$ is {\bf ordinary}, namely if $r=g$, then
$\Fil _0 H^1_{\et}(X_{\overline K},\mu_p)=$
\newline
$H^1_{\et}(X_{\overline K},\mu_p)$,
$\Fil _1 H^1_{\et}(X_{\overline K},\mu_p)=\Ker (\Sp)\simeq 
H^1_{\et}(X_{\overline k},\Bbb Z/p\Bbb Z)$, and we have a canonical
$G_K$-equivariant exact sequence:
$$0\to H^1_{\et}(X_{\overline k},\Bbb Z/p\Bbb Z)\to
H^1_{\et}(X_{\overline K},\mu_p)\to
H^1_{\et}(X_{\overline k},\mu_p)\to 0.$$

\pop {2}
\par
\noindent
{\bf \gr References}

\pop {.5}
\par
\noindent
[Ep] H.Epp, {\sl Eliminating wild ramification}, Invent. Math. 19,
235-249, (1973).

\pop {.5}
\par
\noindent
[Gr-Ma] B. Green M. and M. Matignon, {\sl Order $p$-automorphisms of the 
open disc of a $p$-adic field}. J. Amer. Soc. 12 (1), 269-303, (1999). 

\pop {.5}
\par
\noindent
[Gr] A. Grothendieck, S\'eminaire de g\'eom\'etrie alg\'ebrique SGA-1, 
Lecture Notes 224, Springer Verlag, (1971).

\pop {.5}
\par
\noindent
[He] Y. Henrio, {\sl Arbres de Hurwitz et automorphismes d'ordre $p$ des 
disques et couronnes $p$-adiques formels}. Th\`ese de doctorat, 
Bordeaux France (1999).

\pop {.5}
\par
\noindent
[Hy] O. Hyodo, {\sl Wild ramification in the imperfect residue field case},
Advanced Studies in Pure Mathematics 12, (1987).

\pop {.5}
\par
\noindent
[Mi] Milne, \'Etale cohomology, Princeton University Press, (1980).

\pop {.5}
\par
\noindent
[Oo-Se-Su], F. Oort, T. Sekiguchi and N. Suwa, {\sl On the deformation of 
Artin-Schreier to Kummer}, Ann. Scient. Ec. Norm. Sup. 22, 345-375, (1989).

\pop {.5}
\par
\noindent
[Ra] M. Raynaud, {\sl Sp\'ecialisation du foncteur de Picard}, 
Pub. IHES 38, 27-76, (1970).

\pop {.5}
\par
\noindent
[Sa] M. Sa\"\i di, {\sl Wild ramification and a vanishing cycles formula}, 
math.AG/0106248.

\pop {.5}
\par
\noindent
[Sa-1] M. Sa\"\i di, {\sl Galois covers of degree $p$: semi-stable reduction 
and Galois action}, 
\newline
math.AG/0106249.

\pop {.5}
\par
\noindent
[Se] J. P. Serre, Local Fields, 3 \'edit., Herman, Paris, (1980).

\pop {2}
Mohamed Sa\"\i di

\pop {1}
\par
Departement of Mathematics
\par
University of Durham
\par
Science Laboratories
\par
South Road
\par
Durham DH1 3LE United Kingdom
\par
saidi\@durham.ac.uk

\enddocument
\end